\documentclass[10pt]{amsart}

\usepackage{amsmath, amssymb, amsthm, latexsym,nicefrac}

\theoremstyle{plain}
\newtheorem{theorem}{Theorem}[section]

\newtheorem{lemma}[theorem]{Lemma}
\newtheorem{corollary}[theorem]{Corollary}

\newtheorem{definition}[theorem]{Definition}

\theoremstyle{remark}

\input{xy}
\xyoption{all}

\def\mn{\par\medskip\noindent}
\def\cU{\mathcal U}
\def\cO{\mathcal O}

\def\dim{{\rm dim}}
\def\Zz{\mathbb Z}
\def\Qz{\mathbb Q}
\def\Cz{\mathbb{C}}
\def\Rz{\mathbb{R}}
\def\id{\mathrm{id}}
\def\Nz{\mathbb{N}}
\def\Proj{{\rm Proj}}
\def\disc{{\rm disc}}
\def\fix{{\#_{\rm 1}}}
\def\rank{{\rm rank}}

\title{Sofic groups and diophantine approximation}
\author{Andreas Thom}
\address{Andreas Thom, Mathematisches Institut der Universit\"at G\"ottingen,
Bunsenstr. 3-5, D-37073 G\"ottingen, Germany}
\email{thom@uni-math.gwdg.de}
\urladdr{http://www.uni-math.gwdg.de/thom}

\begin{document}

\begin{abstract}
We prove the algebraic eigenvalue conjecture of J.Dodziuk, P.Linnell, V.Mathai, T.Schick and S.Yates (see \cite{5authors})
for sofic groups. Moreover, we give restrictions on the spectral measure of elements in the integral group ring. Finally, we 
define integer operators and prove a quantization of the operator norm below $2$.
To the knowledge of the author, there is no group known, which is not sofic.
\end{abstract}

\maketitle

\section{Introduction}

The concept of sofic groups goes back to M.Gromov, see \cite{gromov}, who called these groups
\textit{initially sub-amenable}. Later, the term \textit{sofic} was coined by 
B.Weiss in \cite{weiss}, who studied these groups in connection with dynamical systems.

In \cite{elekszabo, elekszabo2, elekszabo3}, G.Elek and E.Szab\'o have extended the study of sofic groups so
that finally, it is known that the class of sofic groups is closed under sub-groups, direct products, extensions with amenable
quotients, direct and inverse limits and free products. In particular, it contains all residually amenable groups. 
To the knowledge of the author, it is not known whether
there are groups which are not sofic. 
\mn
Throughout, $\Gamma$ denotes a discrete group, $\ell^2\Gamma$ the Hilbert space with basis $\Gamma$ and
$L\Gamma$ the von Neumann algebra of those bounded linear operators on $\ell^2\Gamma$, that commute with the right convolution action of
$\Gamma$ on $\ell^2 \Gamma$. There is a natural inclusion $\lambda\colon \Cz \Gamma \to L\Gamma$ of
the complex group ring as left convolution operators. For the theory of von Neumann algebras, we refer to
\cite{tak3} and the references therein.
\mn 
It has attracted much interest to study the spectral properties of the operators $\lambda(a)$ for $a \in \Zz \Gamma$.
An excellent survey and compendium of most of the known results with proofs 
is the book by W.L\"uck, see \cite{lueck}.
\mn
The final goal of this article is to prove the following theorem, which in a similar form 
has been conjectured by Dodziuk, Linnell, Mathai, Schick and Yates, 
as Conjecture $4.14$ in \cite{5authors}, to hold for all groups:
\begin{theorem}
Let $\Gamma$ be a sofic group and $A \in M_n(\Zz \Gamma)$. Consider the natural action of $M_n(\Zz \Gamma)$ on
the Hilbert space $\ell^2\Gamma^{\oplus n}$ by left convolution. Denote by $\lambda(A)$ the 
corresponding bounded operator.
\begin{enumerate}
\item[(i)] All eigenvalues of the $\lambda(A)$ are algebraic integers. 
\item[(ii)] If $\alpha \in \Cz$ is an eigenvalue of the operator $\lambda(A)$, then
all its Galois conjugates are also eigenvalues and the corresponding 
eigenspaces have the same von Neumann dimension.
\end{enumerate}
\end{theorem}

As we will see, 
it is immediate to deduce from the method of proof that the corresponding statements hold for $\Qz \Gamma$ 
and even for $\cO \Gamma$ and 
$\overline{\Qz}\Gamma$, the group rings over the ring of 
algebraic integers $\cO$ and the ring of algebraic numbers $\overline{\Qz}$.

\begin{theorem}
Let $\Gamma$ be a sofic group and $A \in M_n(\overline{\Qz} \Gamma)$. Consider the convolution operator
$\lambda(A)$ as above.
\begin{enumerate}
\item[(i)] All eigenvalues of $\lambda(A)$ are algebraic numbers.
\item[(ii)] If all entries in the coefficient matrices are in $\cO$, then the eigenvalues are algebraic integers.
\item[(iii)] If all entries in the coefficient matrices are in $\Qz$, then all Galois conjugates of eigenvalues
appear and the corresponding eigenspaces have the same von Neumann dimension.
\end{enumerate}
\end{theorem}

These theorems will be proved as Theorem \ref{main} and \ref{main2} in Section \ref{applic}. 
Similar results have been obtained by the authors of \cite{5authors} for amenable groups and groups
in Linnell's class ${\mathcal C}$, see \cite{5authors}.
Having
proved Theorems \ref{main} and \ref{main2}, we continue
to study the resulting restrictions on the spectral measure of elements in the integral group ring. As an
application of several results on zero-distributions of integer polynomials, we derive a quantization of
the operator norm on self-adjoint elements of the integral group ring of a sofic group below $2$.
\mn
The paper is organized as follows. In Section \ref{sofic} we introduce the class of sofic groups
and some auxiliary algebras. Section \ref{approx} deals with diophantine approximation of complex numbers.
We introduce a measure of integer algebraicity, which we show to detect algebraic integers with their degree. 
Section \ref{applic} contains a proof of the algebraic eigenvalue conjecture \cite[Conj. 4.14]{5authors}.
In the remaining two sections, we study the restrictions on the spectral measure, imposed by the study in
Section \ref{applic} and prove the quantization theorem for so-called integer operators. 
Again, several results on diophantine approximation are recalled and applied.
\section{Sofic groups} \label{sofic}
Elek and Szab\'o have given several distinct characterizations of sofic groups, see \cite{elekszabo2}. 
The one that serves our purposes most is the following.
Here, we denote by $S_n$ the permutation group on $n$ letters.
Let $\sigma \in S_n$ be a permutation, we define $\fix \sigma$ to be the number of fixed points of $\sigma$.

\begin{definition} \label{sofic-group}
A group $\Gamma$ is called \emph{sofic} if for any real number $0<\epsilon<1$ and any finite
subset $F\subseteq \Gamma$ there exists a natural number $n \in \Nz$ and a map
$\phi\colon\Gamma\to S_n$ with the following properties:
\begin{enumerate}
\item[(i)]   $\fix \phi(g)\phi(h)\,\phi(gh)^{-1}\ge(1-\epsilon)n,\quad$ for any two elements $g,h\in F$,
\item[(ii)]   $\phi(e)=e \in S_n$,
\item[(iii)]   $\fix\phi(g)\le\epsilon n,\quad$ for any $e\ne g\in F$.
\end{enumerate}
\end{definition}

We may assume that $F$ is symmetric, i.e. $F^{-1}=F$, and $\phi(g^{-1}) = \phi(g)^{-1}$, for $g \in F$. 
Let us order the pairs $(F,\epsilon)$ of a finite symmetric sub-set $F \subset \Gamma$ and a real number $0 < \epsilon <1$
according to $(F,\epsilon) \leq (F',\epsilon')$ if $F\subset F'$ and $\epsilon \geq \epsilon'$. This ordered set
is clearly directed.
Moreover, if $\Gamma$ is sofic, we pick a map $\phi_{F,\epsilon}\colon \Gamma \to S_{n_{F,\epsilon}}$ 
for each $(F,\epsilon)$, satisfying the conditions in the definition above.

Using the standard representation $\pi\colon S_n \to GL_n(\Cz)$ we can extend the maps $\phi_{F,\epsilon}$
to a linear map
$$\phi\colon \Cz\Gamma \to \prod_{F,\epsilon} M_{n_{F,\epsilon}} \Cz.$$
Note that the image of $\Zz \Gamma$ is contained in $\prod_{F,\epsilon} M_{n_{F,\epsilon}} \Zz$.
This map is not multiplicative, but the failure of being multiplicative is under control as $(F,\varepsilon)$
tends to infinity. Indeed, 
$$\frac {\rank\left(\phi_{F,\epsilon}(ab) - \phi_{F,\epsilon}(a)\phi_{F,\epsilon}(b)\right)}{n_{F,\epsilon}}
\to 0,\quad \forall a,b \in \Cz \Gamma,$$ as $(F,\epsilon) \to \infty$. This is obvious if $a,b \in \Gamma$ and 
follows in general since $\phi_{F,\epsilon}$
is linear. Considering an appropriate quotient of $\prod_{F,\epsilon} M_{n_{F,\epsilon}} \Cz$ allows to
deduce direct finiteness of $\Cz \Gamma$ quite easily, moreover: for this to be true, $\Cz$ can be replaced by any skew-field. 
This is the argument which was used by Elek and Szab\'o to prove this result in \cite{elekszabo3}.
\mn
The failure of $\phi$ being multiplicative can be controlled in different ways. Let $A \in \Cz \Gamma$ be fixed, then the family
$(\phi_{F,\epsilon}(A))$ is uniformly operator-norm bounded by $\|A\|_1$. Hence,
$$\|\phi_{F,\epsilon}(A)\|_2 \leq \sqrt{\rank(\phi_{F,\epsilon}(A))} \cdot \|A\|_1$$ is also a good measure, 
as long as the families are uniformly bounded. 
Let us pick a non-trivial ultra-filter $\omega$ 
on the ordered set $$\{(F,\epsilon)\mid F\subset_{\rm fin} \Gamma, 0<\epsilon <1\},$$ 
and let $\lim_{\omega}$ denote the limit along the ultra-filter. Consider the unital complex algebra
$$\left(\prod_{F,\epsilon}{}^b M_{n_{F,\epsilon}} \Cz\right)/ J_{2}$$ which 
consists of uniformly operator-norm bounded families $(a_{F,\epsilon})$ modulo 
the ideal $J_2$ consisting of those families, for which $\|a_{F,\epsilon}\|^2_2/n_{F,\epsilon} \to_{\omega}0$. It is well-known, see
e.g. \cite{tak3}, that $\left(\prod^b_{F,\epsilon} M_{n_{F,\epsilon}} \Cz\right)/ J_{2}$ is a finite von
Neumann algebra, called the ultra-product of the family of algebras $M_{n_{F,\epsilon}}\Cz$,
with respect to the ultra-filter $\omega$. It comes equipped 
with a faithful and normal trace state $$\tau = \lim_{\omega} \tau_{F,\epsilon}.$$ 
Here, $\tau_{F,\epsilon}$ denote the normalized trace on $M_{n_{F,\epsilon}}\Cz$.

Note that the inclusion $\iota\colon \Cz \Gamma 
\to \left(\prod^b_{F,\epsilon} M_{n_{F,\epsilon}} \Cz\right)/ J_{2}$
is a trace and involution preserving homomorphism and hence has a natural injective, trace preserving and normal 
extension to $L\Gamma$. This shows that sofic groups (or rather their group von Neumann algebras) 
satisfy the Connes Embedding Conjecture, which asserts that every tracial von Neumann algebra
can be embedded into an ultra-product of finite-dimensional algebras. 
Again, the argument is due to Elek and Szab\'o, see \cite{elekszabo}.
\mn
%

Every finite tracial von Neumann algebra $(M,\tau)$ carries a notion of rank, i.e.\,for $a \in M$, we define:
$$[a] = \inf\{\tau(p)| p \in \Proj(M), pa=a\}.$$
There is an associated rank metric $d(a,b) = [a-b]$ which equips $M$ with the structure of a metric ring. 
Usually, i.e.\,if $M$ is infinite dimensional, $(M,\tau)$ is not complete with respect 
to that metric, but it follows that its completion is a unital ring. In fact, Lemma \ref{affi} describes
the completion as a well-known companion. This metric and the associated completion functor 
was used by the author in \cite{thom1}, to give a simplified proof of Gaboriau's Theorem about the 
invariance of $L^2$-Betti numbers under orbit equivalence.
\mn
If $(M,\tau)$ is a finite von Neumann algebra, we consider the densely defined closable operators on $L^2(M,\tau)$, 
that are affiliated with $(M,\tau)$, i.e.\,those operators whose real and imaginary part
has all spectral projections in $M$. It is known, that the set of those operators is a unital ring which we call
$\cU(M,\tau)$.

\begin{lemma} \label{affi}
The completion of a finite tracial von Neumann algebra $(M,\tau)$ with respect to its
rank metric is naturally identified with the algebra of closable densely defined operators affiliated with $(M,\tau)$.
Moreover, a trace-preserving inclusion of finite von Neumann algebras extends to
the algebras of affiliated operators.
\end{lemma}
\begin{proof}
The notion of rank has a natural extension to affiliated operators and the algebra of affiliated operators
is complete with respect to the associated rank metric. Hence, it is sufficient to show that $(M,\tau)$ 
is dense in $\cU(M,\tau)$. For this, it is enough to show that all positive affiliated operators $A \in \cU(M,\tau)$ 
are in the closure of $(M,\tau)$. The operator $A$ has a spectral decomposition with respect to which
one can consider the spectral projections $p_{[0,n]}$, associated to the interval $[0,n]$. 
Clearly, $Ap_{[0,n]} \in M$ and since $p_{[0,n]}$ 
increases to one, we see that $Ap_{[0,n]}$ has limit $A$ in rank metric, as $n \to \infty$.
\end{proof}

For brevity, we denote the algebra of operators affiliated with $(L\Gamma,\tau)$ by $\cU\Gamma$,
and the operators affiliated with $\left(\prod^b_{F,\epsilon} M_{n_{F,\epsilon}} \Cz\right)/ J_{2}$
by $\cU_{\omega}(\Gamma)$, suppressing most of the choices we made to define it.

\section{Approximation of complex numbers} \label{approx}
For $\lambda \in \Rz_{>0}$ and $n\in \Nz$, denote by $P_n(\lambda)$ the set 
of monic polynomials of degree at most $n$ with rational integer coefficients, which have all 
their zeros inside $$B(0,\lambda) = \{\mu \in \Cz\mid |\mu| \leq \lambda\}.$$
We set $P(\lambda) = \cup_{n \geq 1} P_n(\lambda)$.
\begin{lemma} \label{disc}
Let $n \in \Nz$ and $\lambda \in \Rz_{>0}$. The set $P_n(\lambda)$ is finite.
\end{lemma}
\begin{proof}
The $k$-th coefficient is bounded by $\binom{n}{k}\lambda^k$. Since we consider only polynomials
with integer coefficients the result follows.
\end{proof}

Let $\beta \in \Cz$ be a complex number in $B(0,\lambda)$. Let $\epsilon>0$ and $p(t) \in P(\lambda)$ be a polynomial.
We set $Z(\beta,\epsilon,p(t))$ to be the number of zeros (counted with multiplicities) of the polynomial $p(t)$
which lie inside $B(\beta,\epsilon)= \{\mu \in \Cz \mid |\mu - \beta| \leq \epsilon\}$.
We are interested in the asymptotics
of the quantity $Z(\beta,\epsilon,p(t))/\deg(p)$ (i.e.\,the proportion of the zeros that lie in an $\epsilon$-neighborhood 
of $\beta$) as $\epsilon$ tends to zero and $n$ tends to infinity, where $p(t)$ is allowed to 
vary among all polynomials in $P_n(\lambda)$.

More concretely, we consider the following real number:
\[ \eta(\beta,\lambda) = \sup \{\delta\geq0 \mid \forall \epsilon >0\, \exists p(t) \in P(\lambda) \colon
Z(\beta,\epsilon,p(t))/\deg(p) \geq \delta \}\]
Let us define:
$$\eta(\beta) = \sup\{\eta (\beta,\lambda)\mid \lambda>0\}.$$

\begin{lemma}\label{easy}
If $\beta$ is an algebraic integer with minimal polynomial $p(t)$, then $$\eta(\beta) \geq {\deg(p(t))}^{-1}.$$
\end{lemma}
\begin{proof} This is clear, since the roots of $p(t)$ are all different and the proportion of roots near $\beta$
is greater than $\deg(p(t))^{-1}$, independent of $\varepsilon.$
\end{proof}
Our main result concerning the quantity $\eta$ is a converse to the preceding Lemma. This will be proved
as Theorem \ref{algeb}. 
\mn
For the proof we need some preparation. Let $p(t)$ be a monic polynomial. The Mahler measure $M(p)$ is defined to be
$$M(p) = \prod_{\alpha\colon p(\alpha)=0} \max(1,|\alpha|),$$
where multiplicities are taken into account.
The following Theorem is inspired by work of D.Roy and M.Waldschmidt. The proof
follows the idea of the proof of Proposition $10.1$ in \cite{waldroy}, although 
the statement is different. We do not need the full force of the following result, however we want to 
state it as follows, since we think it is itself interesting.

\begin{theorem}\label{est}Let $n,t$ be integers and $2 \leq t \leq n$ and let $\beta \in \Cz$ be a complex
number. For any monic polynomial $p(t) \in \Cz[t]$ of 
degree $n$ with distinct zeros and any choice of $t$ distinct zeros $\alpha_1,\dots,\alpha_t$ of $p(t)$, we have
\[\max_{1 \leq i \leq t} |\alpha_i - \beta| \geq  |\disc(p)|^{-\frac1{t(t-1)}} \left( 2 M(p)^{2/n} \right)^{-\frac{n(n-1)}{t(t-1)}} . \]
\end{theorem}
\begin{proof}  We choose an ordering of the zeros of $p(t)$ which satisfies
$$i\leq j \Rightarrow |\alpha_i - \beta| \leq |\alpha_j - \beta|.$$
The discriminant of $p(t)$ is defined to be:
$$\disc(p) = \prod_{i<j} (\alpha_i - \alpha_j)^2.$$
The following estimates can be checked easily:
$$|\alpha_i - \alpha_j| \leq \left\{ \begin{array}{ll} 2 \max\{1,|\alpha_i|\} \max\{1,|\alpha_j|\} |\alpha_j - \beta| & i < j \leq t\\
2 \max\{1,|\alpha_i|\} \max\{1,|\alpha_j|\} & \mbox{otherwise.} \end{array} \right. $$
Using the definition of the discriminant, we get:
\begin{eqnarray*}
 |\disc(p)| &\leq& 2^{n(n-1)} M(p)^{2(n-1)} \prod_{j=1}^t |\alpha_j - \beta|^{2(j-1)} \\
&\leq& 2^{n(n-1)} M(p)^{2(n-1)} |\alpha_t - \beta|^{t(t-1)} \\ 
\end{eqnarray*}
And hence,
\begin{eqnarray*}
\max_{1 \leq i \leq t} |\alpha_i - \beta| &\geq& 
|\disc(p)|^{-\frac1{t(t-1)}} 2^{-\frac{n(n-1)}{t(t-1)}} M(p)^{-\frac{2(n-1)}{t(t-1)}}\\
&=& |\disc(p)|^{-\frac1{t(t-1)}} \left( 2 M(p)^{2/n} \right)^{-\frac{n(n-1)}{t(t-1)}} .
\end{eqnarray*}
\end{proof}

\begin{corollary}
Let $n,t$ be integers and $2 \leq t \leq n$ and let $\beta \in \Cz$ be a complex
number. For any monic irreducible polynomial $p(t) \in \Zz[t]$ of 
degree $n$ and any choice of $t$ distinct zeros $\alpha_1,\dots,\alpha_t$ of $p(t)$, we have
\[\max_{1 \leq i \leq t} |\alpha_i - \beta| \geq  \left( 2 M(p)^{2/n} \right)^{-\frac{n(n-1)}{t(t-1)}} . \]
\end{corollary}
\begin{proof}
Since $p(t)$ is irreducible over $\Zz[t]$, $\gcd(p,p')=1$, 
hence all zeros are distinct. Now, since $\disc(p) \geq 1$, the claim follows from
Theorem \ref{est}
\end{proof}

\begin{theorem}\label{algeb} Let $\beta \in \Cz$ be a complex number.
If $\eta(\beta) \neq 0$, then $\beta \in \Cz$ is an algebraic integer. Moreover, if $\beta$ is an algebraic integer,
than $\eta(\beta) = \deg(p(t))^{-1}$, where $p(t)$ is the minimal polynomial of $\beta$.
\end{theorem}
\begin{proof}
Suppose $\eta(\beta) = \delta \neq 0$. For any 
given $\delta' < \delta$, there exists $\lambda>0$, such that $\eta(\beta,\lambda) > \delta'$.
Hence, for every $0 < \epsilon < 1$, there exists a monic integer 
polynomial $q(t)$, such that the proportion of zeros of $q(t)$ inside $B(\beta,\epsilon)$ is 
greater than $\delta'$. 
Moreover, there is an upper bound $\lambda$ on the moduli of the zeros of $q(t)$, independent of $\epsilon$.
Without loss of generality, we can assume that $q(t)$ is irreducible. Indeed, every monic polynomial over $\Zz[t]$ splits
uniquely into a product of monic irreducible polynomials. If each of these had a proportion of zeros in
$B(\beta,\epsilon)$ less than $\delta'$ , then clearly the product as well.
\mn
We set $n= \deg(q)$ and choose any integer $t' \in \Nz$ with $\delta' n \leq t' \leq n$. Using $M(q) \leq (1 + \lambda)^n$, the notation and the estimate from Theorem \ref{est}, we get
\[ \epsilon \geq \max_{1 \leq i \leq t'} |\alpha_i - \beta| \geq  \left( 2 (1 + \lambda)^2 \right)^{-\frac{n-1}{{\delta'} ({\delta'} n -1)}}
= \left( 2 (1 + \lambda)^2 \right)^{-\frac{1}{{\delta'}^2}\left( 1 + \frac{1 - {\delta'}}{{\delta'} n -1}\right)}. \]
And hence, taking the logarithm:
\[\log \epsilon \geq - \frac{1}{{\delta'}^2}\left( 1 + \frac{1 - {\delta'}}{{\delta'} n -1}\right) \log \left( 2 (1 + \lambda)^2 \right).\]
Solving for ${\delta'}^2$, we get
\[{\delta'}^2 \leq \left( 1 + \frac{1 - {\delta'}}{{\delta'} n -1}\right) \frac{\log \left( 2 (1 + \lambda)^2 \right)}{|\log \epsilon|}.\]
By assumption, for any $\epsilon$ we find $q(t)$, such that the assertion holds. Taking
$$\epsilon \leq \exp\left(-\frac{2 \log ( 2 (1 + \lambda)^2)}{{\delta'}^2}\right),$$ yields
$1 \leq \frac12(1 + \frac{1 - {\delta'}}{{\delta'} n -1})$ and hence $\deg(q) = n \leq \frac{2-{\delta'}}{\delta'}$. 
This implies that $\beta$ is an algebraic integer since $P_n(\lambda)$ is finite 
by Lemma \ref{disc} and all zeros of polynomials in $P_n(\lambda)$ are
algebraic integers. Moreover, for sufficiently small $\epsilon>0$, $\beta$ has to be a root of $q(t)$ and 
$\delta' \leq \deg(q(t))^{-1} \leq \deg(p(t))^{-1}$, where $p(t)$ is the minimal polynomial of $\beta$. Since, $\delta' < \delta$
was arbitrary, this implies
that $\eta(\beta) = \delta \leq \deg(p(t))^{-1}$. We have
already seen in Lemma \ref{easy}, that $\eta(\beta) \geq \deg(p(t))^{-1}$, so the claim follows.
\end{proof}

\section{Applications to $L^2$-invariants}
\label{applic}

Before we come to the main theorem of this paper, let us recall some results, which we need in the sequel.
\mn
A unital ring $R$ is called von Neumann regular, if for every $a \in R$, there exists $b \in R$, such that
$aba=a$. One of the basic results in the theory of von Neumann algebras 
is the following: the ring of operators affiliated with a 
finite von Neumann algebra is von Neumann regular. In fact, John von Neumann himself invented
this notion (under a different name) to abstract this property of the algebras of affiliated operators
of (what was later called) von Neumann algebras.
The following result is well-known and only repeated for completeness.

\begin{lemma}\label{over} A von Neumann regular ring is inverse closed in every over-ring,
i.e. if an inverse of some element exists in an over-ring, than it is already in the ring itself. \end{lemma}
\begin{proof}
Let $R \subset S$ be an inclusion of rings, let $R$ be von Neumann regular and $a \in R$.
Let $b \in R$ be such that $aba=a$. If $ab=ba=1$, then $b$ is
an inverse and hence has to coincide with any inverse in $S$. If the element $a$ is not invertible in $R$, 
then without loss of generality $ab \neq 1$. 
We get $(ab -1)a =0$. So, $a$ is a zero-divisor in $R$ and hence not invertible in $S$. 
\end{proof}

The following theorems are the main results of this paper. We prove the algebraic eigenvalue conjecture
of Dodziuk, Linnell, Mathai, Schick and Yates from \cite{5authors} for sofic groups.

\begin{theorem}\label{main}
Let $\Gamma$ be a sofic group and $A \in M_n(\Zz \Gamma)$.
\begin{enumerate}
\item[(i)] All eigenvalues of the matrix $A$ (acting on $\ell^2\Gamma^{\oplus n}$) are algebraic integers. 
\item[(ii)] For all $\lambda \in \Cz$, the approximation formula holds:
$$\dim_{L\Gamma} \ker(\lambda \cdot \id -A) = \lim_{(F,\epsilon)\to \infty} 
\frac{\dim_{\Cz} \ker(\lambda \cdot \id -\phi_{F,\epsilon}(A))}{n_{F,\epsilon}}.$$
\item[(iii)] If $\alpha \in \Cz$ is an eigenvalue of the matrix $A$, then
all its Galois conjugates are also eigenvalues and the corresponding eigenspaces have the same von Neumann dimension.
\end{enumerate}
\end{theorem}
\begin{proof}
(i) Let us assume $n=1$. 
In order to show that there are only integer algebraic eigenvalues, let us take $\lambda \not \in \cO$ and show
that $\lambda \cdot \id -A$ is invertible in $\cU\Gamma$. This is sufficient, since 
it yields that $\lambda\cdot \id -A$ is not a zero-divisor in $L\Gamma$ and hence the spectral projection onto its kernel 
has to be zero.
\mn
Since $\cU \Gamma \subset \cU_{\omega}(\Gamma)$ (recall the notation from the end of Section \ref{sofic}),
and $\cU \Gamma$ is inverse closed in any over-ring, by Lemma \ref{over}, it is sufficient to construct an inverse in
$\cU_{\omega}(\Gamma).$

The matrix $\lambda \cdot \id - \phi_{F,\epsilon}(A)$ is invertible for all $(F,\epsilon)$, 
since all eigenvalues of $\phi_{F,\epsilon}(A)$ 
are integer algebraic and $\lambda \not \in \cO$
Indeed, the eigenvalues are the roots of the characteristic polynomial of $\phi_{F,\epsilon}(A)$. Since the matrix
has only rational integer coefficients, the characteristic polynomial is a monic polynomial with rational 
integer coefficients. Hence, its roots are algebraic integers.
\mn
Denote by $p_{C,F,\epsilon}$ the biggest projection, such that 
$(\lambda \cdot \id - \phi_{F,\epsilon}(A))^{-1}p_{C,F,\epsilon}$ is of operator-norm less than $C$.
It is sufficient to show that
$$\lim_{C \to \infty} \lim_{\omega} \tau_{F,\epsilon}(p_{C,F,\epsilon}) =1.$$
Indeed, under this condition,
$\left((\lambda \cdot \id - \phi_{F,\epsilon}(A))^{-1}p_{C,F,\epsilon}\right)
\in \prod^b_{F,\epsilon} M_{n_{F,\epsilon}} \Cz/J_2$ converges in rank
metric to the inverse of the operator $\lambda \cdot \id -A$ in $\cU_\omega(\Gamma)$, as $C$ tends to infinity.
\mn
The convergence of the limit to one is achieved by analyzing how the zeros of the characteristic polynomial of the matrix
$\phi_{F,\epsilon}(A)$ approximate $\lambda$. 
Indeed, convergence of $ \lim_{\omega} \tau_{F,\epsilon}(p_{C,F,\epsilon})$ as
$C \to \infty$ is clear, since the $p_{C,F,\epsilon}$ are increasing, let us 
assume that $\lim_{C \to \infty} \lim_{\omega} \tau_{F,\epsilon}(p_{C,F,\epsilon}) = 1 - \delta$.
We get, that for all $C>0$ and $\delta' < \delta$, there exists a $(F,\epsilon)$, such that the 
proportion of zeros of $\chi_{\phi_{F,\epsilon}(A)}$ in $B(\lambda,C^{-1})$ is greater 
than $\delta'$. Moreover, all zeros of those characteristic polynomials are in modulus smaller 
than $\|A\|_1$, since $\|\phi_{F,\epsilon}(A)\| \leq \|A\|_1$. This 
implies $\eta(\lambda)\geq \delta$ and by Theorem \ref{algeb}, we get that $\lambda$ has to be an algebraic integer. 
\mn
(ii) For $\lambda \not \in \cO$, the claim follows from (i). 
Any algebraic integer in the spectrum, has to be a root of a characteristic 
polynomial $\chi_{\phi_{F,\epsilon}(A)}(t)$ for high $(F,\epsilon)$. 
Indeed, as the proof of Theorem \ref{algeb} shows, the zeros of the characteristic 
polynomials cannot accumulate in a proportion $\delta>0$, unless they coincide with the accumulation point.
The same argument as for (i) shows that
the inverse of $(p_{\lambda,F,\epsilon}+\lambda \cdot \id - \phi_{F,\epsilon}(A))$ 
exists in $\cU_{\omega}(\Gamma)$ and hence
$$\dim_{L\Gamma} \ker(\lambda \cdot \id -A) \leq \lim_{(F,\epsilon)\to \infty} 
\frac{\dim_{\Cz} \ker(\lambda \cdot \id -\phi_{F,\epsilon}(A))}{n_{F,\epsilon}}.$$
The other inequality is obvious, since $(p_{\lambda,F,\epsilon})(\lambda \cdot \id -A) =0$
in $\cU_{\omega}(\Gamma)$.
This proves the claim. Note, that $(p_{\lambda,F,\epsilon})$ represents the projection onto the kernel
of $\lambda \cdot \id - A$.
\mn
(iii) The statement concerning the conjugate algebraic integers is now obvious, since 
every zero of the characteristic polynomial comes with its minimal polynomial as a 
factor of the characteristic polynomial and hence with all its conjugate algebraic integers. So, the
right side of the equation in (ii) is invariant under the action of the Galois group.
\end{proof}

The following refinement follows from the methods.

\begin{theorem} \label{main2}
Let $\Gamma$ be a sofic group and $A \in M_n(\overline{\Qz} \Gamma)$.
\begin{enumerate}
\item[(i)] All eigenvalues of $A$ (acting on $\ell^2\Gamma^{\oplus n}$) are algebraic numbers.
\item[(ii)] If all entries in the coefficient matrices are rational, then all Galois conjugates of eigenvalues
appear and the corresponding eigenspaces have the same von Neumann dimension.
\item[(iii)] If all entries in the coefficient matrices are algebraic integers, then the eigenvalues are algebraic integers.
\item[(iv)] The approximation formula holds: $$\dim_{L\Gamma} \ker(\lambda \cdot \id -A) = \lim_{(F,\epsilon)\to \infty} 
\frac{\dim_{\Cz} \ker(\lambda \cdot \id -\phi_{F,\epsilon}(A))}{n_{F,\epsilon}}.$$
\end{enumerate}
\end{theorem}
\begin{proof}
Statement (ii) follows from Theorem \ref{main} by considering a suitable integer multiple of $A$. Similarly,
(i) follows from (iii). It remains to prove (iii) and (iv). \mn (iii) Consider a finite Galois field extension of $\Qz$ in which all the entries
of $A$ exist and consider the Galois conjugates $\sigma_i(A)$ 
of the matrix $A$. Instead of $A$, we consider the matrix $\widehat{A} = A \oplus \sigma_1(A) \oplus 
\dots \oplus \sigma_n(A)$.
Applying the direct sum of $\phi_{F,\epsilon}$ to $\widehat{A}$, yields matrices with characteristic 
polynomials that have rational integer coefficients. Indeed, these coefficients are algebraic integers and invariant
under the Galois action. Moreover, an eigenvalue of $A$ is also an eigenvalue
of $\widehat{A}$, so that it is sufficient to prove that $\widehat{A}$ has only integer algebraic eigenvalues. 
The proof proceeds as before.
\mn
(iv) Recall, using the notation used in the proof of Theorem \ref{main}, 
that $(p_{\lambda,F,\epsilon})$ represents the projection onto the kernel of
$\widehat{A} = A \oplus \sigma_1(A) \oplus 
\dots \oplus \sigma_n(A)$. Moreover, the projection 
$1 \oplus 0 \oplus \dots \oplus 0$ commutes with all diagonal operators. This implies the claim.
\end{proof}
We end this section by giving another application of the methods we used. However, 
a proof is not given, since the theory is already well-developed and the result can be obtained from 
the Determinant Conjecture (see \cite[Thm.\,13.3]{lueck}), 
which is known to be true for sofic groups, \cite[Thm.\,5]{elekszabo}. 
\mn
Let $\Gamma$ be a discrete group and consider a decreasing sequence $H_{n}$, $n \in \Nz$,
of normal subgroups such that $\cap_{n \in I} H_{n} = \{e\}$. Let $A \in M_n(\Zz \Gamma)$ and consider
the images $A_{n} \in M_n(\Zz (\Gamma/H_{n}))$. The matrix $A$ acts naturally
on $\ell^2(\Gamma)^{\oplus n}$ and similarly do the matrices $A_{n}$ act on $\ell^2(\Gamma/H_{n})$. 
The Approximation Conjecture of W.L\"uck asserts that
$$\dim_{L\Gamma} \ker(A) = \lim_{n \to \infty} \dim_{L(\Gamma/H_{n})} \ker(A_{n}), \quad
\forall A \in M_n(\Zz \Gamma).$$
We point out that our methods are strong enough to give a direct proof of the approximation conjecture
in the following form. 

\begin{theorem} \label{approx2}
Let $\Gamma$ be a discrete group and $H_n$ be a decreasing system
of normal subgroups as above, such that the quotients $\Gamma/H_n$ are sofic. 
The approximation conjecture holds, i.e.\,for all $A \in M_n(\Zz \Gamma)$,
$$\dim_{L\Gamma} \ker(A) = \lim_{n \to \infty} \dim_{L(\Gamma/H_{n})} \ker(A_{n}),$$
where we use the notation from above.
\end{theorem}

As before, the techniques (an in particular (iv) of Theorem \ref{main}) 
do also imply the truth of the approximation conjecture for elements in
$\overline{\Qz}\Gamma$, however deducing it for elements in $\Cz \Gamma$ or even $\ell^1 \Gamma$ 
seems to require new ideas.

\section{The spectral measure of integer operators}

We start with the definition of integer operators in an arbitrary finite von Neumann algebra.

\begin{definition}
Let $(M,\tau)$ be a finite von Neumann algebra with a trace. A normal operator $a \in M$ is called
\emph{integer}, if its spectral measure can be approximated weakly by the spectral measures
of uniformly operator-norm bounded sequence of normal integer matrices.
\end{definition}

\begin{lemma}
Let $\Gamma$ be a sofic group. $\lambda(\Zz\Gamma) \subset (L\Gamma,\tau)$ consists of
integer operators.
\end{lemma}

It is easy to see that the conclusions of Theorem \ref{main} carry over verbatim to
integer operators in any finite von Neumann algebra. Indeed, the sub-algebra, which is
generated by the integer operator embeds into an ultra-powers, precisely as in the situation which
was studied.

We are now interested in the
spectral measure of an integer operator $A^*=A$. In particular, we study the distribution function
$R_{A}(\lambda) = \tau(p_{\lambda})$, where $p_{\lambda}$ is the spectral projection
corresponding to the interval $(-\infty,\lambda)$.

We have seen in Theorem \ref{main}, that the size of possible atoms, i.e.\,jumps of the
distribution function, is the limit
of the corresponding sizes of atoms of integer matrices $\phi_{F,\epsilon}(A)$
that approximate $A$.
\mn
Let $p,q \in \Zz[t]$ be a monic irreducible polynomials. We denote by $n(p,q)$ the highest power
of $p$, that appears as a factor in $q$.
We set $$\delta_p = \lim_{\omega} 
\frac {n(p,\chi_{\phi_{F,\epsilon}(A)})}{n_{F,\epsilon} \cdot \deg(p)}.$$

By, Theorem \ref{main}, $A$ will have an atom at 
any root of $p(t) \in \Zz[t]$ of size $\delta_p$, and hence
$\sum_{p} \delta_p \cdot \deg(p) \leq 1$.
Let us set $$q_{F,\epsilon} = \frac{\chi_{\phi_{F,\epsilon}(A)}}
{\prod_{p\colon \delta_p \neq 0} p^{n(p,\chi_{\phi_{F,\epsilon}(A)})}},$$
that is: we are ignoring all polynomials that contribute to atoms. 
The zero distributions of $q_{F,\epsilon}$ converge in distribution (and hence point-wise)
to the continuous part of the spectral
distribution of $A$, suitably normalized. The following Lemma is well-known.

\begin{lemma}[Poly\'a] \label{monot} 
If a sequence of distributions converges point-wise to a continuous distribution, then
it converges uniformly.
\end{lemma}
\begin{proof}
It is sufficient to show that a point-wise convergent 
sequence of monotone functions $$f_n\colon [0,1] \to [0,1],\quad f_n(0) = 1- f_n(1) = 0,$$ 
with continuous limit $f$ is uniformly convergent.
Note that  $f(0) = 0$, $f(1) = 1$, and $f$ is
increasing. Let $\epsilon > 0$. 
By uniform continuity of $f$, choose $k$ 
such that $f(\nicefrac{i}{k}) - f(\nicefrac{i-1}{k}) < \epsilon/2$ for $1 \leq i \leq k$. 
Next, choose $m$ such that if $n \geq m$, then 
$|f(\nicefrac{i}{k}) - f_n(\nicefrac{i}{k})| < \epsilon/2$ 
for all $1 \leq i \leq k$.
If $x$ is in $[0,1]$, then $x$ is in $[\nicefrac{i-1}{k}, \nicefrac{i}{k}]$ for some $i$, so:
$$f_n(x) \leq f_n(\nicefrac{i}{k}) < f(\nicefrac{i}{k}) + \epsilon/2 < f(\nicefrac{i-1}{k}) + 
\epsilon/2 + \epsilon/2 \leq f(x) + \epsilon.$$
Similarly, $f_n(x) > f(x) - \epsilon$, and hence $|f(x) - f_n(x)| < \epsilon$.
\end{proof}

In our case, Lemma \ref{monot} says that the zero distribution 
functions of the polynomials $q_{F,\epsilon}$ converge uniformly to the continuous part of the 
distribution function of $A$. We formulate this observation as Theorem \ref{app}. 
\mn
Recall, an algebraic integer
is called totally real, if all its Galois conjugates are real as well. If $\alpha$ is a totally real algebraic integer,
we denote by $\mu_{\alpha}$ the sum of the Dirac measures of weight $\deg(\alpha)^{-1}$ at the
conjugates of $\alpha$.

\begin{theorem} \label{app}
Let $(M,\tau)$ be a a finite von Neumann algebra and $A=A^*$ an integer operator in $M$. 
The distribution function of the spectral
measure of $A$ is a sum of 
\begin{enumerate}
\item[(i)] a step function $R^a_A$ with jumps at a full set of conjugates of totally real integer algebraic integers, and
\item[(ii)] a continuous function $R^c_A$ which is the uniform limit of distributions of 
convex combinations of measures $\mu_{\alpha_i}$, where $\alpha_i$ is a totally real algebraic integer.
\end{enumerate} 
\end{theorem}

It would be very interesting to characterize the continuous parts of the spectral distribution function intrinsically.
The only restriction on the modulus of continuity that can be derived from the techniques above is 
$$|R^c_A(\lambda) - R^c_A(\mu)| \leq C\left| \log|\lambda-\mu| \right|^{-\frac{1}2},$$
for some constant $C>0$. We mention this result without proof,
since its significance is not clear yet.
A reasonable conjecture for $\Zz\Gamma$ with $\Gamma$ torsion-free would be, 
that there is no singular part of the distribution
function, i.e. absolute continuity of $R_A$.

\section{Quantization of the norm}

The following quantization theorem, shows that the values of the operator-norm on self-adjoint integer operators
in a finite von Neumann algebra are quantized below $2$. 
\mn
We need some preparation. In $1857$, L.Kronecker proved in \cite{kron}, that  if an algebraic integer $\alpha$
lies together with all its Galois conjugates in $[-2,2]$, then is has to be of the form $2\cos(\pi q)$ for some
rational number $q$. Moreover, if $p$ is a monic integer polynomial with all its roots in $[-2,2]$, then $t^{\deg(p)}p(t + \nicefrac 1t)$
is a product of cyclotomic polynomials. 
This is the source of many quantization phenomena in the interplay between integrality
and analysis.

\begin{theorem} Let $(M,\tau)$ be a finite von Neumann algebra and let $A \in M$ be a non-zero integer operator.
\begin{enumerate}
\item[(i)] We always have $\|A\|\geq 1$. 
\item[(ii)] If $A = A^*$, then either
 $\|A\| \geq 2$, or $\|A\|= 2\cos(\pi/q)$, for some natural number $q \in \Nz$. Moreover, in this case 
the spectral measure is completely atomic and the atoms occur at positions $2\cos(\pi p/q)$, for natural numbers $p \in \Nz$.
\item[(iii)] If $A=A^*$ and $\|A\|=2$, then either there is an atom at $2$, or there is a sequence of atoms at $2\cos(\pi/q_n)$
for some sequence of natural numbers with $q_n \to \infty$, or $A$ has full spectrum in $[-2,2]$ and the continuous
part of the measure is semi-circular distributed.
\end{enumerate}
\end{theorem}
\begin{proof} (i) follows from (ii). Indeed $\|A\|^2=\|A^*A\|$ and $A^*A$ is self-adjoint.

(ii) Suppose that $A=A^*, \|A\| <2$ and choose 
$ \|A\| < 2\cos \phi < 2$. Consider the characteristic 
polynomials $\chi_{\phi_{F,\epsilon}(A)}$
of $\phi_{F,\epsilon}(A)$.
Let $d_{F,\epsilon}/n_{F,\epsilon}$ be the proportion 
of zeros of $\chi_{\phi_{F,\epsilon}(A)}$ which
lie outside $[-2\cos \phi,2 \cos \phi]$. It is clear that $d_{F,\epsilon}/n_{F,\epsilon}$ converges to
zero. Indeed, the ultra-limit of this sequence is bounded above by the spectral measure of the complement of
$[-\|A\|,\|A\|]$, which is zero.
We would like to exclude these zeros, but in order to retain integrality, we then would also 
have to exclude all their Galois conjugates. Let us denote by $d'_{F,\epsilon}$ the 
number of those roots that have a conjugate outside $[-2\cos \phi,2 \cos \phi]$.
We have to show that 
$$\lim_{\omega}d'_{F,\epsilon}/n_{F,\epsilon}=0.$$
Indeed, this implies that we can ignore the whole minimal polynomial (as a factor 
$\chi_{\phi_{F,\epsilon}(A)}$) as soon as one root of it lies
outside $[-2\cos \phi,2 \cos \phi]$.
The strategy is to show that the ratio $d_{F,\epsilon}/d'_{F,\epsilon}$ can be bounded from 
below, independent of $(F,\epsilon)$. Again, this will be achieved using special properties of the distribution of zeros
of integer polynomials.
\mn
Let $\alpha$ be a root of $\chi_{\phi_{F,\epsilon}(A)}$, not in $[-2\cos \phi,2 \cos \phi]$. 
Denote by $p(t)$ its minimal polynomial. Transforming everything according to 
$$q(t) = t^{\deg(p)}p(t+ \nicefrac 1t),$$
the roots of $p$ outside $[-2\cos \phi,2\cos \phi]$ transform into $2$ distinct roots of $q$ in the angular sector
$[-\phi,\phi] \cup [\pi - \phi, \pi + \phi ]$. All other roots of $p$ transform into $2$ roots of $q$ that lie on the unit circle. Moreover,
all roots of $q$ are in modulus smaller that $$\lambda=\frac12\left(\|A\|_1 + \sqrt{\|A\|_1^2-4}\right).$$
It is now clear that Lemma \ref{dub} implies our assertion. 
\mn
Being now able to ignore all roots outside $[-2\cos \phi,2 \cos \phi]$, we can use 
Kronecker's Theorem, and conclude that the eigenvalues of $\phi_{F,\epsilon}(A)$ are restricted 
to numbers of the form $2\cos(\pi q)$, for some finite set of rational numbers $q \in \Qz$ as above. This implies the claim.

(iii) Let $\chi_{\phi_{F,\epsilon}(A)}$ be the characteristic polynomial and consider 
$$\psi_{F,\epsilon(t)}(t) = t^{n_{F,\epsilon}}\cdot \chi_{\phi_{F,\epsilon}(A)}(t+\nicefrac 1t).$$
The first two cases correspond to the situation where the proportion of cyclotomic factors in $\psi_{F,\epsilon}$
of fixed degree is relevant as $(F,\epsilon) \to \infty$.
No non-cyclotomic polynomial can arise in a relevant proportion as $(F,\epsilon) \to \infty$, 
since it would contribute an atom outside the unit circle by Kronecker's Theorem.
\mn
We have to determine the spectral measure which arises 
as a limit of zero distributions of monic polynomials $p_n(t)$, with distinct roots, which have all their zeros either on
the unit circle or in the interval $(0,\lambda)$. Moreover, we can assume that the proportion of zeros $\delta_n$ of $p_n(t)$
that do not lie on the unit circle goes to zero. Indeed, assume it does not tend to zero. 
Theorem \ref{main} implies that those zeros cannot contribute to an atom at $1$, hence there would be
spectral measure of the limiting distribution in $(1,\lambda]$, which is not. We conclude that
$$M(p_n) \leq \exp( \log(\lambda) \delta_n \deg(p)),$$ 
for some zero-sequence $\delta_n$. Theorem $2.1$ of \cite{gran}, which goes back to Y.Bilu \cite{bilu}, 
gives the desired conclusion that the limiting distribution has to be the equilibrium 
distribution on the unit circle and hence the claim follows.
\end{proof}

\begin{lemma} \label{dub}
Let $p \in \Zz[t]$ be a monic irreducible polynomial with all roots in $B(0,\lambda)$ for some $\lambda$. 
Let $\phi \in (0,\pi)$. Assume that all zeros of $p(t)$ lie either on the unit circle
or in the sector $\{\alpha \in \Cz \mid \arg(\alpha) \in [-\phi,\phi] \cup [\pi-\phi,\pi+\phi]\}$. 
Then, the proportion of zeros that lie in the sector is greater than some $\delta(\phi,\lambda)>0$, 
independent of $p(t)$.
\end{lemma}

The first result concerning the angular distribution of zeros of integer polynomials is the celebrated result
of P.Erd\"os and P.Tur\'an in \cite{ertu}.
\begin{theorem}[Erd\"os-Tur\'an] Let $p(t) = \sum_{i} a_i t^i \in \Zz[t]$ be a monic and irreducible polynomial and let 
$- \pi < \alpha < \beta < \pi$. Then, 
\[ \left| \frac{\#\{\alpha \in \Cz\mid p(\alpha) =0, \arg(\alpha) \in [\alpha,\beta] \}}{\deg(p)} - \frac{|\alpha - \beta|}{2\pi} \right| 
\leq 16 \sqrt{\frac{\log L(p)}{\deg(p)}},\]
where $L(p) = |a_0|^{-\nicefrac 12}(\sum_{i} |a_i|).$
\end{theorem}

However, this first result has been improved by many people. In order to prove Lemma \ref{dub} 
we use a more sophisticated estimate which was obtained by A.Dubickas as Corollary $3$ in \cite{dub2}.

\begin{theorem}[Dubickas] \label{dubth}
Let $p(t) \in \Zz[t]$ with $\deg(p) \geq 2$ be an irreducible monic polynomial.
If $$0<\frac{\log M(p)}{\deg(p)} < \frac56,$$
then
$$\left|\frac{\#\{\alpha \in \Cz\mid p(\alpha) =0, \arg(\alpha) \in [\alpha,\beta] \}}{\deg(p)} - \frac{|\alpha - \beta|}{2\pi} \right| 
< 6 \left( \frac{\log M(p)}{\deg(p)} \right)^{\frac 13} \cdot \log\left(\frac {\deg(p)}{\log M(p)}\right).$$
\end{theorem}
\begin{proof}[Proof of Lemma \ref{dub}] First of all, the statement is true for cyclotomic polynomials and by
Kronecker's Theorem, these are the only ones with $M(p)=1$. Let $p(t)$ be an arbitrary irreducible non-cyclotomic
polynomial.
Let $\delta$ be the proportion of zeros that lie in the sector $[-\phi,\phi]$. In the argument we will 
assume that we can choose $p(t)$ such that $\delta$ is arbitrarily small and derive a contradiction. This implies, 
that there is a lower bound as desired.
We have
$$M(p) = \prod_{p(\alpha)=0} \max\{1,|\alpha|\} \leq \lambda^{\delta \deg(p)},$$ and since $0<\log M(p)/\deg(p) \leq
\delta \log \lambda < 5/6$, for $\delta$ small enough, the
assumption of Theorem \ref{dubth} is fulfilled.
Theorem \ref{dubth} yields:
\[\left|\delta - \frac{\phi}{\pi}\right| < 6 \left( \frac{\log M(p)}{\deg(p)} \right)^{\frac 13} \cdot \log\left(\frac {\deg(p)}{\log M(p)}\right). \]
Taking $\delta$ small, $\log M(p)/\deg(p)$ will be close to $0$, so 
that finally, for a sufficiently small choice of $\delta$, the right side is smaller than $\phi/2\pi$. Indeed,
this is true since $$\lim_{n \to \infty} n^{-\frac13} \log(n) =0.$$ Hence $\delta> \phi/\pi$ and we arrive at
a contradiction.
\end{proof}

We end this article with an easy corollary of the quantization theorem, for which we do not know any
elementary proof.

\begin{corollary}
Let $(M,\tau)$ be a finite von Neumann algebra and let $A=A^* \in M$ be an integer operator. Either,
\begin{enumerate}
\item[(i)] $\|A\|=1$ and $A^2$ is a projection, or 
\item[(ii)] $\|A\| \geq \sqrt{2}$.
\end{enumerate}
\end{corollary}

Just to illustrate the power of the methods and the peculiarity of the possible results, let us give 
one more result which we can obtain from a well-known estimate on the distribution of
zeros of integer polynomials. 
\mn
Recall that for a finite von Neumann algebra $M$, we set:
$$\det{}_M(A) = \exp\left(\frac12 \int_{0^+}^{\infty} \ln(x) \,d\mu_{A^*A}(x)\right),$$
where $\mu_{A^*A}$ denote the spectral measure of $A^*A$. This quantity is called the Fuglede-Kadison
determinant of $A$. For information about it, we refer to \cite[Ch. 3.2]{lueck}.

\begin{theorem} Let $\Gamma$ be a sofic group and let $A=A^* \in M_n\Zz \Gamma$, such that
$0, \pm 1$ are not atoms of spectral measure of $A$. We consider the operator $t \cdot \id + A \in M_n \Zz[\Zz \times \Gamma]$.
We have,
$$\det{}_{L(\Zz \times \Gamma)} (t \cdot \id + A) \geq \sqrt {\frac{1+ \sqrt{5}}2 } = 1,2720 \dots .$$
\end{theorem}
\begin{proof}
This follows from the 
Schick's proof of Theorem $6.9$ in \cite{schick1} and a result of C.Smyth in \cite{smyth}.

Indeed, Schick showed, as an application of Fatou's Lemma, that the Fuglede-Kadison 
determinant is bounded below by the 
limes superior of the absolute values of the determinants $\det_{L[\Zz]}(t \cdot \id - \phi_{F,\epsilon}(A))$ 
of the the matrices, that approximate $A$ in distribution. Since $0,\pm1$ are not atoms of $A$, we can ignore 
by Theorem \ref{main} those atoms in the spectral measure of $\phi_{F,\epsilon}(A)$. Smyth's Theorem says that
the Mahler measure of the characteristic polynomial of $\phi_{F,\epsilon}(A)$, being a polynomial
with only real integer algebraic roots different from $-1,0,1$ is greater than the square-root of 
$\frac{1 + \sqrt5}2$. Now, L\"uck has shown in Example $3.22$ of \cite{lueck} that the $L\Zz$-determinant of a matrix
$t\cdot \id -A$ over $\Zz[\Zz]$ is precisely the Mahler measure of the characteristic polynomial of $A$.
This finishes the proof.
\end{proof}
\mn
\begin{center} \textsc{Acknowledgment} \end{center}
I want to thank Thomas Schick and Ulrich Bunke for constant encouragement and interesting discussions.

\end{document}